\documentclass[12pt,reqno]{amsart}

   \usepackage{amsthm}
   \usepackage{amsmath}
   \usepackage{amssymb}
   \usepackage{amsfonts}
   \usepackage{amscd}
   \usepackage[mathscr]{eucal}
   \usepackage{verbatim}

\newtheorem{thm}{Theorem}[section]

\newtheorem{prop}[thm]{Proposition}
\newtheorem{cor}[thm]{Corollary}
\newtheorem{rem}[thm]{Remark}
\newtheorem{defi}[thm]{Definition}

\newcommand{\adots}{\mathinner{\mkern1mu\raise1pt\hbox{.}
\mkern2mu\raise4pt\hbox{.}
\mkern2mu\raise7pt\hbox{.}\mkern1mu}}

\newcommand{\G}{{\bf G}}

\renewcommand{\H}{{\bf H}}

\newcommand{\T}{{\bf T}}

\newcommand{\Ad}{\mathbb A}
\newcommand{\cx}{\mathbb C}
\newcommand{\rt}{\mathbb Q}

\newcommand{\rl}{\mathbb R}

\renewcommand{\aa}{\alpha}

\newcommand{\bb}{\beta}

\newcommand{\g}{\gamma}

\newcommand{\gl}{\ensuremath{\mbox{GL}}}
\newcommand{\gso}{\ensuremath{\mbox{GSO}}}
\newcommand{\gsp}{\ensuremath{\mbox{GSp}}}
\newcommand{\gspin}{\ensuremath{\mbox{GSpin}}}

\renewcommand{\l}{\lambda}

\newcommand{\oo}{\omega}

\newcommand{\s}{\sigma}

\newcommand{\w}[1]{\ensuremath{\widetilde{#1}}}

\begin{document}

\title{Generic Transfer from $\gsp(4)$ to $\gl(4)$}
\author[Mahdi Asgari]{Mahdi Asgari$^{^\star}$}
\address{School of Mathematics \\
Institute for Advanced Study \\
1 Einstein Drive \\
Princeton, NJ 08450 \\
USA}
\email{asgari@math.ias.edu}
\thanks{$^\star$ Partially supported by the NSF grant DMS--0111298 at IAS}
	
\author[Freydoon Shahidi]{Freydoon Shahidi$^{^\dag}$}
\address{Mathematics Department \\ 
Purdue University \\
West Lafayette, IN 47907 \\
USA}
\email{shahidi@math.purdue.edu}
\thanks{$^\dag$ Partially supported by the NSF grant DMS--0200325}

\begin{abstract}
We establish Langlands functoriality for the generic 
spectrum of $\gsp(4)$ and describe its transfer on $\gl(4)$. 
We apply this to prove results toward the generalized 
Ramanujan conjecture for generic representations of $\gsp(4)$. 
\end{abstract}

\maketitle

\section{Introduction}\label{intro}

Let $k$ be a number field and let $G$ denote the group $\gsp(4,\Ad_k)$. 
The (connected component of the) $L$-group of $G$ is $\gsp(4,\cx)$ 
which has a natural embedding into $\gl(4,\cx)$. Langlands functoriality 
predicts that associated to this embedding there should be a transfer 
of automorphic representations of $G$ to those of $\gl(4,\Ad_k)$ \cite{arthur-gsp4}. 
Langlands' theory of Eisenstein series reduces the proof of this to 
unitary cuspidal automorphic representations. We establish functoriality 
for the generic spectrum of $\gsp(4,\Ad_k)$. More precisely, 
(cf. Theorem \ref{main}) we prove: 

{\it Let $\pi$ be a unitary cuspidal representation of $\gsp(4,\Ad_k)$ 
which we assume to be globally generic. Then $\pi$ has a unique transfer 
to an automorphic representation $\Pi$ of $\gl(4,\Ad_k)$. The transfer 
is generic (globally and locally) and satisfies $\oo_\Pi = \oo_\pi^2$ and 
$\Pi \simeq \w{\Pi}\otimes\oo_{\pi}$. Here, $\oo_\pi$ and $\oo_\Pi$ denote 
the central characters of $\pi$ and $\Pi$, respectively.}

Moreover, we give a {\it cuspidality criterion} for $\Pi$ and prove that, 
when $\Pi$ is not cuspidal, it is an isobaric sum of two unitary cuspidal 
representations of $\gl(2,\Ad_k)$ (cf. Proposition \ref{cusp-cri}). 

We give a number of applications of this result. The first one is 
Theorem \ref{est} which gives estimates toward the generalized 
Ramanujan conjecture for generic 
representations of the group $\gsp(4)$ (cf. Section \ref{ram-est}). We also 
prove in Theorem \ref{weak} that any generic unitary cuspidal representation 
of $\gsp(4,\Ad_k)$ is {\it weakly Ramanujan} (cf. Section \ref{ram-weak} for 
definition). In Section \ref{spinor} we use our main result to give an 
immediate proof that the spin $L$-function of a generic unitary cuspidal 
representation of $\gsp(4,\Ad)$ is entire. This is due to the fact that 
the spin $L$-function of $\gsp(4)$ now becomes a standard $L$-function 
(or product of two such $L$-functions) of general linear groups. 
This fact has also been proved recently by R. Takloo-Bighash using different 
methods from ours. 

We should note here that the transfer from $\gsp(4)$ to $\gl(4)$ has been 
expected by experts in the field for a long time. It is our understanding 
that Jacquet, Piatetski-Shapiro, and Shalika knew how to prove this 
result, at least in principle, but, as far as we know, their result was 
never published. We should also point out that their proof is based on 
methods that are fairly disjoint from ours.  

Our method of proof is to start with our earlier, more general but weaker, 
result on generic transfer from $\gspin$ groups to $\gl$ (cf. \cite{gspin}). 
This gives us the existence of $\Pi$. We then use results of Piatetski-Shapiro and 
Soudry on analytic properties of $L$-functions of $\gsp(4)$ twisted by 
$\gl(1)$ and $\gl(2)$ to get more information about the representation 
$\Pi$. It is exactly the lack of such results in the general case of 
$\gspin$ groups that prevents us from carrying out our analysis for the 
more general case for now. However, as we pointed out in \cite{gspin}, 
there are currently two ways to overcome this problem. One is to prove 
an analogue of descent theory for these groups as was done for classical 
groups by Ginzburg, Rallis, and Soudry \cite{grs,sou}. The other 
is to use techniques along the lines of \cite{kim1,kim2}. 

The case of the transfer for all automorphic representations of 
$\gsp(4)$ (whether generic or not) requires Arthur's trace formula, 
whose process in this case and what the expected issues are, is 
outlined in \cite{arthur-gsp4}.
 
We would like to thank James Cogdell for many helpful discussions and 
Peter Sarnak for his interest in this work. The first author would 
also like to thank Steve Gelbart, Dihua Jiang, Robert Langlands, Brooks 
Roberts, and Ramin Takloo-Bighash for many helpful discussions during 
the course of this work.

\section{Main Result} 

Let $k$ be a number field and let $\Ad = \Ad_k$ denote its ring of adeles.  
We define the similitude symplectic group of degree four via  
\[ \gsp(4) = \left\{g \in \gl(4) \, : \, ^tg J g = \mu(g) J \right\}, \]
where 
\[ J= \left(\begin{matrix}  
&&& 1 \\ && 1 & \\ & -1 && \\ -1 &&& \end{matrix}  \right) \]
and $\mu(g)\in\gl(1)$ is the similitude character. We fix the 
following parametrization of the elements of the maximal torus 
$\T$ in $\gsp(4)$: 
\[ \T = \left\{ t = t(a_0,a_1,a_2) = 
\left(\begin{matrix} 
a_0 a_1 a_2 &&& \\ & a_0 a_1 && \\ && a_0 a_2 & \\ &&& a_0 
\end{matrix} \right) 
\right\}. \]
The above agrees with our previous more general notation for the 
group $\gspin(2n+1)$ in \cite{gspin}. Recall that 
the group $\gsp(4)$ is identified with $\gspin(5)$.

Let $\pi = \otimes^\prime_v \pi_v$ be a globally $\psi$-generic unitary cuspidal 
automorphic representation of $\gsp(4,\Ad)$. Here, $\psi = \otimes_v \psi_v$ 
is a non-trivial additive character of $k\backslash\Ad$ defining a character of 
the unipotent radical of the standard upper-triangular 
Borel in the usual way. 
We fix $\psi$ now and do not repeat it in the rest 
of this paper. Let $S$ be any non-empty 
finite set of non-archimedean places $v$ which includes those 
$v$ with $\pi_v$ or $\psi_v$ ramified. We proved in \cite{gspin} 
that there exists an automorphic representation 
$\Pi = \otimes^\prime_v \Pi_v$ of $\gl(4,\Ad)$ such that $\Pi_v$ is 
a local transfer of $\pi_v$ for $v$ outside of $S$. 

To be more explicit, assume that $v\not\in S$. If $v$ is archimedean, then $\pi_v$ 
is given by a parameter $\phi_v : W_v \longrightarrow \gsp(4,\cx)$ 
where $W_v$ is the Weil group of $k_v$ (cf. \cite{langlands-real}).
Let $\Phi_v : W_v \longrightarrow \gl(4,\Ad)$ be given by 
$\Phi_v = \iota \circ  \phi_v$, where $\iota : \gsp(4,\cx) \longrightarrow \gl(4,\cx)$ 
is the natural embedding. Then $\Phi_v$ is the parameter of $\Pi_v$. 

If $v \not\in S$ is non-archimedean, then $\pi_v$ is the unique unramified 
subquotient of the representation induced from an unramified character 
$\chi$ of $\T(k_v)$ to $\gsp(4,k_v)$. Writing 
\begin{equation}\label{unram} 
\chi(t(a_0,a_1,a_2)) = \chi_0(a_0) \chi_1(a_1) \chi_2(a_2), 
\end{equation}
where $\chi_i$ are unramified characters of $k_v^\times$ and $a_i \in k_v^\times$ 
the representation $\Pi_v$ is then the unique irreducible unramified 
subquotient of the representation of $\gl(4,k_v)$ parabolically induced 
from the character 
\begin{equation}\label{chi} 
\chi_1\otimes\chi_2\otimes\chi_2^{-1}\chi_0\otimes\chi_1^{-1}\chi_0 
\end{equation}
of $\T(k_v)$.

Moreover, we proved that $\oo_\Pi = \oo^2$, where $\oo = \oo_\pi$ 
and $\oo_\Pi$ denote the central characters of $\pi$ and $\Pi$, respectively, 
and for $v\not\in S$ we have $\Pi_v\simeq\w{\Pi}_v\otimes\oo_{\pi_v}$, i.e., 
$\Pi$ is nearly equivalent to $\w{\Pi}\otimes\oo$.  

The representation $\Pi$ is equivalent to a subquotient of some representation 
\begin{equation}\label{aut-ind} 
\mbox{Ind} (|\det|^{r_1}\s_1 \otimes \cdots \otimes |\det|^{r_t}\s_t)
\end{equation}
where induction is from $\gl(n_1)\times\cdots\times\gl(n_t)$ with 
$n_1+\cdots+n_t = 4$ to $\gl(4)$ and $\s_i$ are unitary cuspidal automorphic 
representation of $\gl(n_i,\Ad)$ and $r_i\in\rl$. 

Without loss of generality we may assume that $r_1\ge r_2\ge \cdots \ge r_t$. 
Moreover, since $\Pi$ is unitary we have $n_1 r_1 + \cdots + n_t r_t = 0$ 
which implies that $r_t \le 0$. Let $T = S \cup \{v : v|\infty\}$ and 
consider 

\begin{equation}\label{partial}
L^T(s,\pi \times \w{\s}_t) = L^T(s,\Pi\times \w{\s}_t) = 
\prod_{i=1}^t L^T(s+r_i, \s_i \times \w{\s}_t). 
\end{equation}

If $n_t=1$, then the left hand side is entire by a result of Piatetski-Shapiro 
(cf. page 274 of \cite{ps-gsp4}). Now consider the right hand side at 
$s_0=1-r_t\ge 1$. The last term in the product has a pole at $s_0$ while all the 
others are non-zero there since $\Re(s_0+r_i) = 1 + r_i - r_t\ge 1$. 
This is a contradiction. 

Now assume that $n_t=3$, i.e., $t=2$ with $n_1=1$ and $n_2=3$. Replacing 
$\pi$ and $\Pi$ by their contragredients will change $r_i$ to $-r_i$ and 
takes us back to the above situation which gives a contradiction again. 

Therefore, $n_t = 2$. In this case, the left hand side of (\ref{partial}) may 
have a pole at $s=1$ (cf. Theorem 1.3 of \cite{ps-soudry} and beginning of its 
proof), and if so, arguing as above, we conclude that $r_t=0$. This means that 
we either have $t=2$ with $n_1 = n_2 = 2$ or $t=3$ with $n_1 = n_2 = 1$ and $n_2 = 2$. 
However, we can rule out the latter as follows. 

Assume that $t=3$ with $n_1=n_2=1$ and $n_2=2$. Then, it follows from the 
fact that $r_3 = 0$ 
and the conditions $r_1 \ge r_2 \ge r_3$ and $r_1 + r_2 + 2 r_3 = 0$ that 
all the $r_i$ would be zero in this case. This implies that if we consider 
the $L$-function of $\pi$ twisted by $\w{\s}_1$, we have 
\begin{equation}\label{two-blocks}
L^T(s,\pi \times \w{\s}_1) = L^T(s, \s_1 \times \w{\s}_1) 
L^T(s, \s_2 \times \w{\s}_1) L^T(s, \s_3 \times \w{\s}_1).
\end{equation} 
Now the left hand side is again entire by Piatetski-Shapiro's result 
mentioned above and the right hand side has a pole at $s=1$ which is a contradiction. 

Therefore, the only possibilities are $t=1$ (i.e., $\Pi$ unitary cuspidal) or 
$t=2$ and $n_1 = n_2 = 2$ with $r_2 = 0$. In the latter case we immediately 
get $r_1 = 0$ as well since $r_1 + r_2 = 0$ by 
unitarity of the central character. Moreover, in this case 
we have $\s_1 \not\simeq \s_2$ since, otherwise, 
\begin{equation}
L^T(s,\pi \times \w{\s}_1) = L^T(s, \s_1 \times \w{\s}_1) L^T(s, \s_2 \times \w{\s}_1) 
\end{equation}
must have a double pole at $s=1$ while any possible pole of the left hand side 
at $s=1$ is simple (cf. proof of Theorem 1.3 of \cite{ps-soudry}). 

Therefore, we have proved the following: 

\begin{prop}\label{ni2}
Let $\pi$ be a globally generic unitary cuspidal automorphic representation 
of $\gsp(4,\Ad)$ and let $\Pi$ be any transfer of $\pi$ to $\gl(4,\Ad)$. Then, 
$\Pi$ is a subquotient of an automorphic representation as in (\ref{aut-ind}) 
with either $t=1$, $n_1 = 4$, and $r_1 = 0$ 
(i.e., $\Pi$ is unitary cuspidal) or $t=2$, $n_1 = n_2 = 2$ and 
$r_1 = r_2 = 0$. In the latter case, we have $\s_1\not\simeq\s_2$.
\end{prop}

In fact, we can get more precise information.

\begin{prop}\label{cusp-cri} 
Let $\pi$ be a globally generic unitary cuspidal automorphic 
representation of $\gsp(4,\Ad)$ with $\oo=\oo_\pi$ its central character 
and let $\Pi$ be any transfer as above. Then, $\Pi \simeq \w{\Pi}\otimes\oo$ 
(not just nearly equivalent). 
\begin{itemize}
\item[(a)] The representation $\Pi$ is cuspidal if and only if $\pi$ is not obtained 
as a Weil lifting from $\gso(4,\Ad)$. 
\item[(b)] If $\Pi$ is not cuspidal, then it is the isobaric sum of two 
representations $\Pi = \Pi_1 \boxplus \Pi_2$, where each $\Pi_i$ 
is a unitary cuspidal automorphic representation of $\gl(2,\Ad)$ 
satisfying $\Pi_i \simeq \w{\Pi}_i \otimes \oo$ and 
$\Pi_1\not\simeq\Pi_2$. 
\end{itemize}
\end{prop}

\begin{proof}
By Proposition \ref{ni2}, $\Pi$ is not cuspidal if and only if it is a subquotient of 
\begin{equation}\label{Sigma}
 \Sigma=\mbox{Ind}_{\gl(2,\Ad)\times\gl(2,\Ad)}^{\gl(4,\Ad)} 
(\s_1 \otimes \s_2), 
\end{equation}
where $\s_i$ are unitary cuspidal automorphic representation of $\gl(2,\Ad)$. 

On the other hand, by Theorem 1.3 of \cite{ps-soudry} mentioned above the representation $\pi$ 
is obtained as a Weil lifting from $\gso(4,\Ad)$ if and only if there exists 
an automorphic representation $\tau$ of $\gl(2,\Ad)$ such that $L^T(s,\pi\times\tau)$ 
has a pole and in that case $\tau$ can be normalized so that the pole occurs at $s=1$. 

Now assume that $\Pi$ is cuspidal. Then, for any $\tau$ as above we have 
\[ L^T(s,\pi\times\tau) = L^T(s,\Pi\times\tau) \]
which is entire. Therefore, $\pi$ is not obtained as a Weil lifting from $\gso(4,\Ad)$. 
Moreover, since $\Pi$ is cuspidal, so is $\w{\Pi}\otimes\oo$ and they are 
nearly equivalent, therefore, by strong multiplicity one theorem 
\cite{jac-sha-classificationI,jac-sha-classificationII,ps-corvallis} we 
have $\Pi \simeq \w{\Pi}\otimes\oo$.

Next, assume that $\Pi$ is not cuspidal and, hence, is given as a subquotient 
of $\Sigma$ above. We claim that each $\Sigma_v = \mbox{Ind} (\s_{1,v}\otimes\s_{2,v})$ 
is irreducible. To see this note that each $\s_{i,v}$ in generic unitary and 
is either a tempered representation of $\gl(2,k_v)$ or a complementary series 
$I(\chi |\ |^\aa,\chi |\ |^{-\aa})$ with $\aa\in(0,1/2)$ and $\chi$ a 
unitary character. If both of the $\s_{i,v}$'s are tempered, then irreducibility 
of $\Sigma_v$ is clear. If both are complementary series of the form 
$I(\chi_1 |\ |^\aa,\chi_1 |\ |^{-\aa})$ and 
$I(\chi_2 |\ |^\bb,\chi_2 |\ |^{-\bb})$ with $\aa,\bb\in(0,1/2)$ and $\chi_i$ 
unitary characters, then for $\Sigma_v$ to be reducible we should have 
$\aa\pm\bb=\pm 1$ which is not possible. Finally, if one of $\s_{i,v}$'s is 
tempered and the other is complementary series, then we either have 
$\Sigma_v = \mbox{Ind} (\chi |\ |^{-\aa},\chi_1,\chi_2,\chi |\ |^\aa)$,  
$\Sigma_v = \mbox{Ind} (\chi |\ |^{-\aa},
Q(\chi_1 |\ |^{-1/2},\chi_1 |\ |^{1/2}),\chi |\ |^\aa)$, or 
$\Sigma_v = \mbox{Ind} (\chi |\ |^{-\aa},\eta,\chi |\ |^\aa)$. Here, 
$\chi_1,\chi$ are unitary characters, $\aa\in(0,1/2)$, $\eta$ is 
a unitary supercuspidal representation of $\gl(2,k_v)$, and 
$Q(\chi_1 |\ |^{-1/2},\chi_1 |\ |^{1/2})$ denotes the Steinberg 
representation twisted by the unitary character $\chi_1$. Again, in 
all these cases the representation $\Sigma_v$ is irreducible 
since $\aa\in(0,1/2)$. 

Therefore, at every place $v$ the representation $\Pi_v$ is the same as 
the irreducible $\Sigma_v$ and, hence, it is the Langlands quotient. 
This means that $\Pi$ is an isobaric 
representation, i.e., the isobaric sum of $\s_1$ and $\s_2$. 
Again by strong multiplicity one theorem, which remains valid for isobaric 
representations \cite{jac-ps-shalika}, we conclude that $\Pi \simeq
\w{\Pi}\otimes\oo$. Now, just 
take $\Pi_i$ to be $\s_i$. 

Finally, by Proposition 7.4 of \cite{gspin}, which was based on classification 
theorems of Jacquet and Shalika, we know that we either have 
$\Pi_i\simeq\w{\Pi}_i\otimes\oo$ for $i=1,2$ or we have 
$\Pi_1\simeq\w{\Pi}_2\otimes\oo$ (or equivalently, 
$\Pi_2\simeq\w{\Pi}_1\otimes\oo$). However, the latter case will not 
occur when $\pi$ is cuspidal and generic since, otherwise, $\pi$ will 
be nearly equivalent to an Eisenstein series representation, i.e., $\pi$ 
will be a CAP representation of $\gsp(4,\Ad)$. This is impossible by 
Theorem 1.1 of \cite{ps-soudry-five}. This completes the proof. 
\end{proof}

\begin{rem}\label{unique}
Notice that any other transfer $\Pi^\prime$ of $\pi$ is also a subquotient of 
$\Sigma$ in (\ref{Sigma}) which is irreducible. 
Therefore, $\pi$ has a {\it unique} transfer to $\gl(4,\Ad)$ which we continue 
to denote by $\Pi$. In particular, this implies that we have not lost any information 
at the places where we did not have a natural candidate for the local transfer.

Moreover, since $\Pi$ is either a unitary cuspidal representation 
of a general linear group or an 
isobaric sum of two such, every local representation $\Pi_v$ is 
full induced and generic. 

Furthermore, if $\Pi$ is not cuspidal, then $\Pi_1\not\simeq\Pi_2$ 
implies that $L(s,\Pi_1\times\w{\Pi}_2)$ has no pole at $s=1$. 
This implies that the Fourier coefficient 
of $\Pi$ along the unipotent radical of our fixed Borel is non-vanishing, 
i.e., $\Pi$ is globally generic \cite{shahidi81AMJ}. 
\end{rem}

We collect the above results in the following theorem which is our main result. 

\begin{thm}\label{main}
Let $\pi$ be a globally generic unitary cuspidal automorphic representation of 
$\gsp(4,\Ad)$ with central character $\oo$. Then $\pi$ has a unique transfer 
$\Pi$ to $\gl(4,\Ad)$ which satisfies $\Pi \simeq \w{\Pi}\otimes\oo$ and 
its central character is $\oo^2$. Moreover, $\Pi$ is either unitary cuspidal 
or an isobaric sum $\Pi_1\boxplus\Pi_2$ of two inequivalent 
unitary cuspidal automorphic representations 
of $\gl(2,\Ad)$ satisfying $\Pi_i\simeq\w{\Pi}_i\otimes\oo$. The latter is 
the case if and only if $\pi$ is obtained as a Weil lifting from $\gso(4,\Ad)$. 
Furthermore, $\Pi$ is globally generic, i.e., it has a non-vanishing Fourier 
coefficient along the unipotent radical of our fixed Borel subgroup. 
\end{thm}

In fact, we can get more information about the local 
representations at places $v\in S$.

\begin{prop}\label{exp} 
Fix $v\in S$ and let 
\begin{equation}
\pi_v \simeq \mbox{Ind}(\pi_{1,v}|\det|^{b_{1,v}}\otimes \cdots 
\otimes\pi_{t,v}|\det|^{b_{t,v}} \otimes \pi_{0,v}) 
\end{equation} 
be an irreducible generic representation of $\gsp(4,k_v)$, where 
each $\pi_{i,v}$ is a tempered representation of some $\gl(n_i,k_v)$, 
$b_{1,v} > \cdots > b_{t,v}$, and $\pi_{0,v}$  is a tempered generic 
representation of some $\gsp(2m,k_v)$. (Note that there are very 
few possibilities for $n_i$'s and $m$ since $n_1+\cdots n_t + m = 2$. 
We are allowing $m=0$ and, by convention, $\gsp(0) = \gl(1)$. ) 
Let $\oo_v$ denote the central character of $\pi_v$. 
Assume that $\pi_v$ is the local component of the globally generic 
unitary cuspidal representation $\pi$ of $\gsp(4,\Ad)$ and let 
$\Pi$ be its transfer to $\gl(4,\Ad)$. Then the local component 
$\Pi_v$ of $\Pi$ at $v$ is generic and of the form 
\begin{equation}\label{Piv}
\begin{matrix}
\Pi_v & \simeq & \mbox{Ind}\, \Big(\pi_{1,v}|\det|^{b_{1,v}}\otimes \cdots 
\otimes \pi_{t,v}|\det|^{b_{t,v}} \otimes 
\Pi_{0,v} \otimes 
\\
&& (\w{\pi}_{t,v}\otimes\oo_v)|\det|^{-b_{t,v}} \otimes \cdots \otimes 
(\w{\pi}_{1,v}\otimes\oo_v)|\det|^{-b_{1,v}} 
\Big), 
\end{matrix} 
\end{equation}
where, $\Pi_{0,v}$ is a tempered generic representation of $\gl(2m,k_v)$ 
if $m > 0$. 
\end{prop}

\begin{proof} 
Let us remark that, as in Section 7 of \cite{ckpss-classical}, one 
could define the notion of local transfer and obtain complete information 
about such transfers for a general irreducible admissible generic 
representation, whether a local component of a global representation 
or not. (In fact, the representation $\Pi_{0,v}$ would then be the 
local transfer of $\pi_{0,v}$.) However, we do not need the full 
extent of such results in this paper. 

Recall that we already proved (cf. Remark \ref{unique}) that 
each $\Pi_v$ is generic and is full 
induced. 

Let $v\in S$ and consider $\pi_v$ and $\Pi_v$ as in the proposition. 
We first show that if $\rho_v$ is any supercuspidal representation of 
$\gl(r,k_v)$, then 
\begin{equation}\label{gamma}
\g(s,\pi_v\times\rho_v,\psi_v) = \g(s,\Pi_v\times\rho_v,\psi_v).
\end{equation} 
The key here is the fact that there exists a unitary cuspidal 
representation $\rho$ of $\gl(r,\Ad)$ such that its local component 
at $v$ is $\rho_v$ and at all other finite places $w\not=v$ the local 
component $\rho_w$ is unramified (cf. Proposition 5.1 of \cite{shahidi:90annals}). 
Now applying converse theorem of Cogdell and Piatetski-Shapiro with 
$S^\prime = S - \{v\}$ will give the result exactly as in the proof 
of Proposition 7.2 of \cite{ckpss-classical}. Moreover, by multiplicativity 
of $\g$-factors, we conclude that (\ref{gamma}) also holds if 
$\rho_v$ is a discrete series representation of $\gl(r,k_v)$. 

Next, we claim that if $\pi_v$ is tempered, then so is $\Pi_v$. Here, 
again the main tool is multiplicativity of $\g$-factors and the proof 
is exactly as in Lemma 7.1 of \cite{ckpss-classical}. This proves the 
proposition for the case $m=2$. If $m=1$, then the group $\gsp(2m)=\gsp(2)$ 
is the same as $\gl(2)$ and we set $\Pi_{0,v}$ to be $\pi_{0,v}$ itself. 
For $m=0$ we need no choice of $\Pi_{0,v}$. Now, let 
$T=\{w_0\}$ consist of a single finite place $w_0\not=v$ 
at which $\pi_v$ is unramified and consider the 
representation $\Pi^\prime$ of $\gl(4,\Ad)$ whose 
local components are the same as $\Pi$ outside of 
$S$ and are the irreducible induced representations on 
the right hand side of (\ref{Piv}) when $v\in S$. We 
can now apply converse theorem again to $\Pi^\prime$ 
and $T=\{w_0\}$ to conclude that $\Pi^\prime$ is a 
transfer of $\pi$. The key here is that the induced representations 
on the right hand side of (\ref{Piv}) have the right 
$L$-functions. Therefore, by uniqueness of the transfer 
we proved earlier we have $\Pi^\prime_v\simeq\Pi_v$ 
for $v\in S$. This completes the proof.  
\end{proof}

\section{Applications}

We first recall that the current formulation of the Ramanujan 
conjecture for generic cuspidal representations states that 
for any quasi-split group $\H$ and any 
globally generic unitary cuspidal automorphic representation 
$\pi = \otimes_v^\prime \pi_v$ the local components $\pi_v$ are 
tempered for all places $v$. As an application of our main 
theorem we can prove two types of results in this direction: 
estimates toward this conjecture for the group $\gsp(4)$ as well 
as a weaker version of it for generic representations of this group. 

\subsection{Estimates toward Ramanujan}\label{ram-est}

Following \cite{ckpss-classical} we introduce 
the following notation in order to prove estimates. 
Let $\Pi = \otimes_v^\prime \Pi_v$ be a unitary cuspidal 
automorphic representation of $\gl(m,\Ad_k)$. For every place 
$v$ the representation $\Pi_v$ is unitary generic and can 
be written as a full induced representation 
\begin{equation}
\Pi_v \simeq \mbox{Ind}(\Pi_{1,v}|\det|^{a_{1,v}}\otimes \cdots 
\otimes\Pi_{t,v}|\det|^{a_{t,v}})
\end{equation} 
with $a_{1,v} > \cdots > a_{t,v}$ and each $\Pi_{i,v}$ 
tempered \cite{vogan,Zel}. 

\begin{defi} 
We say $\Pi$ satisfies $H(\theta_m)$ with $\theta_m \ge 0$ if 
for all places $v$ we have $-\theta_m \le a_{i,v} \le \theta_m$. 
\end{defi}

The classification of generic unitary dual of $\gl(m)$ \cite{tadic,vogan} 
trivially gives $H(\frac 1 2 )$. The best result currently known 
for a general number field $k$ says that any unitary cuspidal 
representation of $\gl(m,\Ad)$ satisfies 
$H(\frac 1 2 - \frac 1 {m^2+1})$ \cite{luo-rudnick-sarnak}. When 
$k=\rt$ and $m\le 4$ it is $H(\frac 1 2 - \frac 1 {1+ m(m+1)/2})$. 
The same bound is also available for $m > 4$ for $k=\rt$  provided 
that one knows that the symmetric square $L$-function of $\Pi$ 
is absolutely convergent for $\Re(s) > 1$ (cf. \cite{kim-sarnak}) 
but this is only available presently for $m \le 4$. When 
$m=2$ we have the better bounds of $H(1/9)$ for a general 
number field $k$ \cite{kim-shahidiDuke} and $H(7/64)$ 
for $k=\rt$ \cite{kim-sarnak}. 

The Ramanujan conjecture demands $H(0)$. 

Similarly, if $\pi = \otimes_v^\prime \pi_v$ is a unitary 
generic cuspidal automorphic representation of $\gsp(2n,\Ad_k)$, 
then by \cite{muic} and \cite{vogan} each $\pi_v$ can be written 
as a full induced representation 
\begin{equation}
\pi_v \simeq \mbox{Ind}(\pi_{1,v}|\det|^{b_{1,v}}\otimes \cdots 
\otimes\pi_{t,v}|\det|^{b_{t,v}} \otimes \tau_v), 
\end{equation} 
where each $\pi_{i,v}$ is a tempered representation of some 
$\gl(n_i,k_v)$ and $\tau_v$ is a tempered generic representation 
of some $\gsp(2m,k_v)$ with $n_1+\cdots + n_t + m = n$. 

\begin{defi}
We say 
$\pi$ satisfies $H(\theta_n)$ with $\theta_n \ge 0$ if for 
all places $v$ we have $-\theta_n \le b_{i,v} \le \theta_n$. 
\end{defi}

The classification of generic unitary dual of $\gsp(4)$ 
(cf. \cite{lapid-tadic-muic}, for example) trivially 
gives the estimate $H(1)$. 
The Ramanujan conjecture demands $H(0)$ again. 
For a survey of results in this direction and their applications 
we refer to \cite{sarnak-fields-notes,shahidi-borel}.

\begin{thm}\label{est}
Let $k$ be a number field and 
assume that all unitary cuspidal representations of 
$\gl(4,\Ad_k)$ (respectively, $\gl(2,\Ad_k)$) satisfy $H(\theta_4)$ 
(respectively, $H(\theta_2)$) and $\theta_2 \le \theta_4$. 
Then, any globally generic unitary cuspidal representation 
$\pi$ of $\gsp(4,\Ad_k)$ satisfies $H(\theta_4)$. 

If $\pi$ 
transfers to non-cuspidal representation of $\gl(4,\Ad_k)$ 
(cf. Theorem \ref{main}), then it satisfies the possibly better 
bound $H(\theta_2)$.  
\end{thm} 

\begin{proof}
Let $\Pi$ be the functorial transfer of $\pi$ to $\gl(4,\Ad_k)$. 

If $v$ is an archimedean place of $k$, then this is clear since 
in this case local functoriality is well understood through Langlands 
parametrization (cf. proof of Theorem 6.1 of \cite{gspin}, 
for example, for more details). 

Let $v$ be a non-archimedean place of $k$ at which $\pi_v$ is 
unramified. Then, it follows from (\ref{chi}) that 
$\pi_v$ is given by its 
Frobenius-Hecke (Satake) parameter 
which is of the form 
\begin{equation}\label{Pi_vParameter}
\mbox{diag} (\chi_1(\varpi), \chi_2(\varpi), 
\chi_2^{-1}(\varpi)\chi_0(\varpi), \chi_1^{-1}(\varpi)\chi_0(\varpi)), 
\end{equation}
where $\varpi$ denotes a uniformizer of $k_v$. 
If $\Pi$ is cuspidal, then for $i=1,2$ we have by assumption  
\[ q_v^{-\theta_4} \le |\chi_i(\varpi)| \le q_v^{\theta_4}. \] 
If $\Pi$ is not cuspidal, then we have similar inequalities 
with $\theta_4$ replaced by even better estimate of $\theta_2$. 
Since $\pi_v$ is unitary, we have $|\chi_0(\varpi)| = 1$. Therefore, 
Frobenius-Hecke parameters of $\pi_v$ also satisfy 
similar inequalities. 

Next, assume that $v$ is a place of $k$ in $S$. Then, by 
Proposition \ref{exp}, a similar argument as above works again. 
\end{proof}

\begin{cor}
Let $\pi$ be a globally generic unitary cuspidal representation 
of $\gsp(4,\Ad_k)$. Then $\pi$ satisfies $H(15/34)$. If 
$\pi$ transfers to a non-cuspidal representation of $\gl(4,\Ad_k)$, 
then it satisfies $H(1/9)$. If $k=\rt$, then we have the 
better estimates of $H(9/22)$ and $H(7/64)$, respectively.
\end{cor}

\begin{proof}
The proof is immediate if we combine Theorem \ref{est} with 
the known results on estimates for general linear groups mentioned 
above. 
\end{proof}

\begin{cor}
The Ramanujan conjecture for unitary cuspidal representations 
of $\gl(4)$ and $\gl(2)$ imply the Ramanujan conjecture for 
the generic spectrum of $\gsp(4)$.
\end{cor}

\subsection{Weak Ramanujan}\label{ram-weak}

Following \cite{cog-ps-unitarity,ramakrishnan-mrl,kim-ex2} 
we recall the following definition. 

\begin{defi} 
Let $\G$ be a split reductive group over the number field $k$. 
Let $\pi = \otimes^\prime_v \pi_v$ be an automorphic representation 
of $\G(\Ad_k)$. We say that $\pi$ is weakly Ramanujan if given 
$\epsilon > 0$ there exists a set $T$ of places of $k$ containing 
the archimedean ones and the non-archimedean ones with $\pi_v$ 
ramified such that $T$ has density zero and for $v\not\in T$ 
the Frobenius-Hecke parameter $\mbox{diag}(\l_{v,i})$ 
of $\pi_v$ satisfies 
\[
\max_i \{ |\l_{v,i}|, |\l_{v,i}^{-1}| \} \le q_v^\epsilon. 
\]
Here, $q_v$ denotes the cardinality of the residue field. 
\end{defi} 

We will be concerned with the cases of $\G=\gl(m)$ or $\G=\gsp(4)$ 
in this paper. We recall that (unitary) cuspidal representations 
of $\gl(m)$ for $m \le 4$ are weakly Ramanujan (cf. \cite{ramakrishnan-mrl} 
and Propositions 3.7 and 6.3 of \cite{kim-ex2}). 

Let $\pi$ be a globally generic unitary cuspidal representation 
of $\gsp(4,\Ad_k)$. For any $v\not\in T$ as above, let 
\begin{equation} 
\mbox{diag}(a_{0,v},a_{1,v},a_{2,v})
\end{equation} 
be the Frobenius-Hecke parameter 
of $\pi_v$ (cf.(\ref{unram})). Then, as in (\ref{Pi_vParameter}), 
the parameter of the local transfer $\Pi_v$ is given by 
\begin{equation}
\mbox{diag}(a_{1,v},a_{2,v},a_{2,v}^{-1} a_{0,v},a_{1,v}^{-1} a_{0,v}). 
\end{equation}
Moreover, $|a_{0,v}| = 1$ since $\pi_v$ is unitary. Therefore, the above 
results about weak Ramanujan property of unitary cuspidal representations 
of $\gl(m)$ immediately imply the following. 

\begin{thm}\label{weak}
Let $\pi$ be a globally generic unitary cuspidal representation of 
$\gsp(4,\Ad_k)$. Then $\pi$ is weakly Ramanujan.
\end{thm}

\subsection{Spinor $L$-function for $\gsp(4)$}\label{spinor} 
As another application we get the following immediate corollary of our 
main result, Theorem \ref{main}. 

\begin{prop} 
Let $\pi$ be a globally generic unitary cuspidal representation 
of $\gsp(4,\Ad_k)$. Then the spinor $L$-function $L(s,\pi, spin)$ is entire. 
\end{prop}

\begin{proof} 
Let $\Pi$ be the transfer of $\pi$ to $\gl(4,\Ad)$. If $\Pi$ is unitary 
cuspidal, then $L(s,\pi, spin) = L(s,\Pi)$ and if $\Pi=\Pi_1\boxplus\Pi_2$ 
is the isobaric sum of two unitary cuspidal representations of $\gl(2,\Ad)$, 
then we have $L(s,\pi, spin) = L(s,\Pi_1) L(s,\Pi_2)$. In either case the 
$L$-functions on the right hand side are standard representations of the 
general linear group and are entire. 
\end{proof}

\begin{rem} This result has also been proved by R. Takloo-Bighash 
in \cite{ramin} among other things. His methods are different from ours 
and are based on integral representations. 
\end{rem}

\bibliography{mahdi}
\bibliographystyle{plain}

\end{document}